\documentclass[a4paper,12pt,leqno]{amsart} 
\usepackage{amsmath,amsthm,amssymb}  
\usepackage{mathdots}
\usepackage{graphicx}

\numberwithin{equation}{section}

\theoremstyle{plain}
\newtheorem{theorem}{Theorem}[section]
         
\newtheorem{corollary}[theorem]{Corollary}

\theoremstyle{definition}  
\newtheorem{example}[theorem]{Example} 
\newtheorem{remark}[theorem]{Remark} 
\newcommand{\C}{\mathbb C}   
\newcommand{\R}{\mathbb R}
\newcommand{\Z}{\mathbb Z}

\newcommand{\al}{\alpha}
\newcommand{\be}{\beta} 
\newcommand{\ga}{\gamma}
\newcommand{\de}{\delta}
\newcommand{\la}{\lambda}
\newcommand{\si}{\sigma} 

\newcommand{\eps}{\epsilon}
\newcommand{\Om}{\Omega}
\newcommand{\De}{\Delta}
\renewcommand{\th}{\theta}
\newcommand{\om}{\omega}
\newcommand{\Th}{\Theta}
\newcommand{\Ga}{\Gamma}
\newcommand{\ze}{\zeta}

\DeclareMathOperator{\diag}{diag}

\DeclareMathOperator{\End}{End}
\DeclareMathOperator{\Aut}{Aut}
\DeclareMathOperator{\Ker}{Ker}

\DeclareMathOperator{\Ad}{Ad}
\DeclareMathOperator{\ad}{ad}

\newcommand{\SL}{\textrm{SL}}
\renewcommand{\sl}{\frak s\frak l}
\newcommand{\SO}{\textrm{SO}}

\newcommand{\SU}{\textrm{SU}}

\newcommand{\g}{{\frak g}}
\newcommand{\h}{{\frak h}}

\renewcommand{\c}{{\frak c}}

\newcommand{\no}{\noindent}
      
\newcommand{\pr}{\prime} 
\newcommand{\prr}{{\prime\prime}} 
\newcommand{\st}{\ \vert\ }   
\renewcommand{\ll}{\lq\lq}
\newcommand{\rr}{\rq\rq\ }
\newcommand{\rrr}{\rq\rq}  
\renewcommand{\b}{\partial}

\newcommand{\bp}{\begin{pmatrix}} 
\newcommand{\ep}{\end{pmatrix}} 
\newcommand{\bsp}{\left(\begin{smallmatrix}} 
\newcommand{\esp}{\end{smallmatrix}\right)}

\newcommand{\zbar}{  {\bar z}  }

\newcommand{\tbar}{  {\bar t}  }
\newcommand{\ttb}{ {t\bar t}  }

\renewcommand{\i}{ {\scriptscriptstyle\sqrt{-1}}\, }
\newcommand{\ii}{ {\scriptstyle\sqrt{-1}}\, }

\newcommand{\cN}{\mathcal{N}}
\newcommand{\cS}{\mathcal{S}}

\newcommand{\hd}{\h^\pr}
\newcommand{\Ded}{\De^\pr}
\newcommand{\bed}{\be^\pr}

\begin{document}     

\title[Topological-antitopological fusion]{Topological-antitopological fusion and the quantum cohomology of Grassmannians
}  
   
\author{Martin A. Guest}    

\date{}   

\maketitle 

\begin{abstract}We suggest an explanation for the part of the Satake Correspondence which relates the quantum cohomology of complex Grassmannians and the quantum cohomology of complex projective space, as well as their respective Stokes data, based on the original physics approach using the tt* equations. We also use the Stokes data of the tt* equations to provide a Lie-theoretic link between particles in affine Toda models and solitons in certain sigma-models. Along the way, we illustrate some (well known) relations between the tt* equations, the non-abelian Hodge Correspondence, and quantum cohomology.
\end{abstract}

\section{Introduction}\label{intro}

An article \cite{Bo95a} by Bourdeau (with a similar title),  which appeared in the physics literature some 25 years ago, presented a relation between complex Grassmannians and complex projective spaces based on supersymmetric field theory. The topological-antitopological fusion equations (tt* equations), which had just been introduced by Cecotti and Vafa (\cite{CeVa91},\cite{CeVa92}), provided the means to compare these two kinds of spaces: for each of them there is a sigma-model (an example of a  $2$-dimensional $N=2$ supersymmetric field theory), and hence a ground state metric (a solution of the tt* equations), and then physical considerations implied a relation between the metric for the Grassmannian and the metric for complex projective space.  The argument was sketched first 
by Cecotti and Vafa in section 8.3 and Appendix A of \cite{CeVa92}, then the results were made more explicit in \cite{Bo95a}.

As far as we know, the tt* aspects of this result were never taken up in the mathematical literature. However, it is a consequence of the above observation that the quantum cohomology of the Grassmannian can be expressed as the exterior product (in an appropriate sense) of  the quantum cohomology of complex projective space.  This fact --- considered as a quantum version of the geometric Satake Correspondence --- has been studied by mathematicians from the viewpoint of representation theory, algebraic geometry, and mirror symmetry. We refer to section 7 of the article \cite{CDGXX} by Cotti-Dubrovin-Guzzetti for a recent survey\footnote{We also recommend  \cite{CDGXX} for an introduction to 
the geometric Satake Correspondence in the wider sense, which is a categorical relation between intersection cohomology of Schubert varieties in the affine Grassmannian of a Lie group $G$ and representations of the Langlands dual group $G^\vee$.  A precise statement of this relation can be found in section 7 of
\cite{CDGXX}, with references to the original articles.
}
of this interesting work containing clear statements.
It is also known that the monodromy/Stokes data of the respective quantum differential equations are related in a similar way, i.e.\ by taking the exterior product. However, in the mathematical literature on the 
Grassmannian/projective space
correspondence, the \ll real structure\rr implicit in the tt* metric has, so far, not played any role.  

It is our aim to explain how {\em the tt* aspects are in fact the most fundamental.} This will serve two purposes: to give a mathematical foundation for the results from physics (as in \cite{Bo95a}), and to give a conceptual explanation for the mathematical results (as in \cite{CDGXX}). 
What makes this possible is the recent series of papers
\cite{GIL1},\cite{GIL2},\cite{GIL3},\cite{MoXX} which
establish the existence of the relevant solutions of the tt* equations and some of their remarkable properties. The  tt* origin of the correspondence becomes visible when these results are interpreted Lie-theoretically, as in \cite{GH1},\cite{GH2}.   

A benefit of keeping track of the tt* metric is that it represents additional analytic data beyond the purely algebro-geometric data of quantum cohomology. Rather than expecting an equivalence of quantum cohomology, one should expect an equivalence of quantum cohomology together with this extra data (as in the original physics approach). In particular this can reduce the ambiguities which bedevil a purely algebraic approach.

Our explanation may also shed light on the relation between the affine Toda model and various sigma-models.  
It was observed by Freeman and others (see 
\cite{Fr91},\cite{Do91},\cite{Do92} and the review \cite{Co99}) that rather intricate properties of roots of Lie algebras play a significant role in the physics of affine Toda theory, at the quantum level as well as the classical level.  
On the other hand Lerche and collaborators (\cite{LeWa91},\cite{FLMW91}) studied \ll polytopic models\rrr, a type of sigma-model in which the polytope spanned by the weights of a representation contains essential physical data such as vacuum/ground states, particles, mass and spin. Both situations illustrate Zamolodchikov's fundamental idea that this kind of data emerges purely on the grounds of symmetry --- the \ll conformal bootstrap\rrr. 
We shall show that both structures arise from (and are therefore linked by) the Stokes data of the tt* equations in a natural way.

This article is organised as follows.
Section \ref{tt*} presents a broad brush picture of the relevant tt* equations and their physical meaning, just to the extent needed to explain their relevance to the Satake Correspondence.  Section \ref{results} contains the main mathematical results. 
Based on this, section \ref{physics} makes precise the statements in 
section \ref{tt*}, with further details of the physical and mathematical significance of the tt* equations, and a new treatment of particles in Toda models and polytopic models. This is illustrated by the case of the Grassmannian sigma-model, our main example.  

It may be appropriate to make some preliminary remarks on {\em why} the tt* equations govern quantum cohomology, and why the Grassmannian case illustrates this so well.   For this it is important to have in mind the physical background, i.e.\ the fact that solutions of the tt* equations actually represent deformations of conformal field theories,
a process which will be described in sections \ref{morephysics} (physically) and \ref{morett*} (mathematically).  The field theory corresponding to the quantum cohomology of a K\"ahler manifold is generally believed to admit such a deformation, but amongst all solutions of the tt* equations these are isolated examples.  The Grassmannians are, of course, very special K\"ahler manifolds.  Thus, the theory of the tt* equations is not specifically about Grassmannians, or even K\"ahler manifolds, but about more general objects.  These objects are in some sense \ll holomorphic\rr (or, in physical language, topological). 

One general principle which relates holomorphic objects to solutions of the tt* equations is the Hitchin-Kobayashi Correspondence, or non-abelian Hodge Correspondence, 
which relates Higgs bundles to flat connections via harmonic bundles.  How this is relevant to the tt* equations is explained in section \ref{morett*}. Essentially, the tt* metric provides the canonical harmonic representative.  From this viewpoint Grassmannians are related to projective spaces because the same solution of the tt* equations represents both.

Another general principle concerns the relation between solutions of the tt* equations and their asymptotics.
Physics suggests that the deformation of a conformal field theory corresponding to a solution of the tt* equations is uniquely determined by its initial point (the conformal field theory, e.g.\ quantum cohomology), or, equivalently, by its final point (another field theory).  Needless to say, this is neither a general phenomenon for nonlinear differential equations, nor for the Hitchin-Kobayashi Correspondence.  It is related to a homogeneity or \ll phase-invariance\rr condition of the type which occurs in the theory of Frobenius manifolds.

A key idea of Dubrovin \cite{Du93}, relevant to both of these general principles, was to formulate the tt* equations as an isomonodromy condition for a family of flat meromorphic connections.  This makes the tt* equations susceptible to
the Riemann-Hilbert method (which converts the nonlinear tt* p.d.e. to a linear integral equation --- see \cite{FIKN06}). By its very nature the Riemann-Hilbert method has the advantage of relating the solutions directly to their asymptotics (which p.d.e.\ methods, including the Hitchin-Kobayashi Correspondence, generally do not).  Thus the Riemann-Hilbert point of view is close in spirit to the original physics problem, and the asymptotic data can be expected to contain the essential physical or geometrical information.

On the other hand, it is technically difficult to implement the Riemann-Hilbert method. In the very special case of the Grassmannian, and its {\em small} quantum cohomology, the tt* equations become a version of the two-dimensional Toda equations. Here (as we explain in this article) Lie-theoretic techniques are extremely helpful.  Moreover, the Grassmannians are examples of minuscule flag manifolds, and for such manifolds the
cohomology/representation dictionary is particularly simple (see
\cite{CDGXX},
\cite{GoMaXX}, \cite{KaXX}).
It can be expected that future refinements of these methods will lead to results for the quantum cohomology of other K\"ahler manifolds,
and for more general solutions of the tt* equations.

{\em Acknowledgements:}  
The author was partially supported by JSPS grant 18H03668.

\section{The tt*-Toda equations and their physical interpretation}\label{tt*}

Cecotti and Vafa gave several concrete examples of tt* equations, 
one of the main examples being a special case of the affine (i.e.\ periodic) Toda equations, which we call the tt*-Toda equations. Although the tt*-Toda equations exist for any complex Lie algebra, we focus on those of type $A_n$, where complete mathematical results are available. 

\subsection{The tt*-Toda equations of type $A_n$}

The equations 
(formula (7.4) in \cite{CeVa91})
are 
\begin{equation}\label{ttt}
 2(w_i)_{\ttb}=-e^{2(w_{i+1}-w_{i})} + e^{2(w_{i}-w_{i-1})}, \ 
 i\in\Z
\end{equation}
for functions $w_i:\C^\ast\to\R$, where, for all $i$,  

\no (a) $w_i=w_{i+n+1}$ (periodicity),  

\no(b) $w_i=w_i(\vert t\vert)$
(radial condition), and 

\no(c) $w_i+w_{n-i}=0$
(anti-symmetry).  

It turns out (see Theorem \ref{thm1}) that the solutions can be parametrized by points of a compact region in a finite-dimensional vector space.  The asymptotic behaviour of a solution at $t=0$ and at $t=\infty$ can be expressed explicitly in terms of the corresponding  parameter value (Theorems  \ref{thm1} and \ref{thm3}). 
These theorems constitute the mathematical foundation for what follows. It should be emphasized that they apply to the tt*-Toda equations (and may eventually extend to all tt* equations), but not to the Toda equations in general; they are not routine consequences of differential equation theory. Undoubtedly the physical interpretation of the tt*-Toda equations plays a significant role here.

\subsection{Physical interpretation}

In physics, solutions of (\ref{ttt}) describe massive deformations of two-dimensional $N=2$ supersymmetric conformal field theories. For a general introduction to this we recommend \cite{Za96}. Here we just mention the basic ingredients relevant to (\ref{ttt}), postponing more details to sections 
\ref{morephysics} and \ref{morett*}. 

The function $w=(w_0,\dots,w_n)$ represents a (positive-definite) Hermitian metric on a (trivial) vector bundle of rank $n+1$ over the space $\C^\ast$ of \ll couplings\rrr, called the ground state metric, or tt* metric. The fibre of the vector bundle is the space of \ll vacua\rrr. 

The particular form of the equations (\ref{ttt}), together with the periodicity property, is a consequence of having tt* equations of {\em Toda type,} but (b) and (c) are general features of tt* geometry. Namely, the radial property (b) is a homogeneity property (\ll chiral charge conservation\rrr), and 
the anti-symmetry property (c) is equivalent to preservation of a \ll topological metric\rr on the vector bundle.

For a given solution of (\ref{ttt}), the point $t=0$ (the conformal point, or ultra-violet limit) represents the superconformal theory that is being deformed.  This theory is massless.  For each nonzero $t$ the corresponding theory has mass, hence the terminology.  The point $t=\infty$ (the infra-red limit) does not represent a superconformal theory; its meaning will be discussed later.  

Thus, the picture described so far is of a family of deformations of superconformal theories,  each represented by a solution of equation (\ref{ttt}).  The particular physical characteristics of the deformation are revealed through its asymptotics at $t=0$ and at $t=\infty$. As we shall explain later, these are the \ll chiral charges\rr (at $t=0$) and the \ll soliton multiplicities\rr (at $t=\infty$). 

\subsection{Functoriality and the Satake Correspondence}

It is well known (see, for example, \cite{Ba15}) that the two-dimensional Toda equations can be defined for any complex simple Lie algebra $\g$, and that these equations are integrable in the sense that they admit a zero-curvature representation. The corresponding tt*-Toda equations were formulated in \cite{GH2}, as follows.

Let $\h$ be a Cartan subalgebra of $\g$, with root system $\De$, and let $l=\dim_\C\h$.
Let
$
\g = \h \oplus\left(  \oplus_{\be\in\Delta} \ \g_\be \right)
$
be the root space decomposition.
Let $B$ be a positive scalar multiple of the Killing form such that $B(e_\be,e_{-\be})=1$ for root vectors $e_\be$, $\be\in\Delta$.
We define $H_\be\in\h$ by 
$B(x,H_\be)=\be(x)$ for all $x\in\h$.
Thus, for any choice of simple roots $\Pi=\{\be_1,\dots,\be_l\}$ we obtain a basis $H_{\be_1},\dots,H_{\be_l}$ of $\h$.  
The tt*-Toda equations apply to functions $w=w(t,\bar t)$ which take values in the vector space $\h_\sharp=\oplus_{i=1}^l\R H_{\be_i}\cong \R^l$.  

To write the zero-curvature form of the equations, we define
\[
E_-=\sum_{i=0}^l  \tfrac{1}{\sqrt{q_i}} e_{-\be_i},\quad
E_+=\sum_{i=0}^l  \sqrt{q_i} e_{\be_i},
\]
where $\be_0=-\psi$ and  $\psi=\sum_{i=1}^l q_i \be_i$ is the
highest root. We define $q_0=1$.  (In the case $\g=\sl_{n+1}\C$ all $q_i$ are $1$.) 
Then we define $\al=\al^\prime dt + \al^{\prime\prime}d\bar t$ by
\begin{equation}\label{galpha}
\al^\prime = w_t + \tfrac1\lambda \Ad(e^w) E_-,
\quad
\al^{\prime\prime}= -w_{\bar t} + \lambda \Ad(e^{-w}) E_+,
\end{equation}
$\lambda\in\mathbb C^\ast$ being a parameter.
The zero-curvature condition $d\al+\al\wedge\al=0$ is equivalent to
$2w_{t \bar t} =[\Ad(e^w)E_-,\Ad(e^{-w})E_+]$, which reduces to
\begin{equation}\label{gttt}
\textstyle
2w_{t \bar t} = -\sum_{i=0}^l  e^{ -2\be_i(w)} H_{\be_i}.
\end{equation}
This is just a real form of the affine Toda equations; the tt*-Toda equations require additional conditions
$w=w(\vert t\vert)$ and $\si(w)=w$, 
where $\si:\g\to\g$ is a certain involution, first defined by Hitchin in \cite{Hi92} (see  \cite{GH2}, and also \cite{AF95},\cite{Ba15}).
In the case $\g=\sl_{n+1}\C$ we have 
$\sigma(X)=-\Delta \, X^T\,  \Delta$
where $\De$ is the
anti-diagonal matrix $\De=(\de_{i,n-i})_{0\le i\le n}$.

To write (\ref{galpha}) using matrices
in the case $\g=\sl_{n+1}\C$, let us take the standard diagonal Cartan subalgebra with diagonal entries $x_0,\dots,x_n$, and simple roots $\be_1=x_0-x_1,\dots,\be_n=x_{n-1}-x_n$. 
We define $\be_0=x_n-x_0$ (i.e.\ minus the highest root). As root vector $e_{x_i-x_j}$ we take the matrix $E_{i,j}$ which has $1$ in the $(i,j)$ entry and $0$ elsewhere.
Then we have:
\begin{equation}\label{malpha}
 \al = (w_t+\tfrac1\lambda W^T)dt + (-w_{\bar t}+\lambda W)d\bar t,
\end{equation}
where
\[
 w=\diag(w_0,\dots,w_n),\ \
 W=
 \left(
\begin{array}{c|c|c|c}
\vphantom{(w_0)_{(w_0)}^{(w_0)}}
 & \!e^{w_1\!-\!w_0}\! & &
 \\
\hline
  &  & \  \ddots\   & \\
\hline
\vphantom{\ddots}
  & &  &  e^{w_n\!-\!w_{n\!-\!1}}\!\!\!
\\
\hline
\vphantom{(w_0)_{(w_0)}^{(w_0)}}
\!\! e^{w_0\!-\!w_n} \!\!  & &  &  \!
\end{array}
\right).
\]
The zero-curvature condition is equivalent to
$2w_{t\bar t} = [W^T,W]$.  Thus we recover equation (\ref{ttt}), and, if
$w_i=w_i(\vert t\vert)$ and $w_i+w_{n-i}=0$, conditions (a),(b),(c) are satisfied.  

We can observe from this excursion into Lie theory that the connection form of the tt*-Toda equation arises in two steps:  

{\em Step 1} is the choice of the Lie algebra $\g=\sl_{n+1}\C$, giving (\ref{galpha}),
then 

{\em Step 2} is the choice of the standard matrix representation $\la_{n+1}$ of $\g$ on the vector space $\C^{n+1}$,  giving the matrix version (\ref{malpha}).

\no Obviously any other faithful representation $\theta$ of $\g$ gives an equivalent system, with exactly the same solutions. However, {\em  the flat connection (and its geometrical/physical interpretation) depends on $\theta$.}
We claim that the relation between Grassmannians and projective spaces
arises from choosing $\theta=\wedge^k \la_{n+1}$ in step 2 and varying $k$. In other words, from the {\em same} solution of (\ref{ttt}) we obtain {\em different} geometrical interpretations.

In the context of the Satake Correspondence, which in general relates \ll cohomology\rr and \ll representations\rrr, one could say that a solution of (\ref{ttt}) produces a package of data (including quantum cohomology, for certain special solutions) which behaves functorially with respect to a choice of representation.

We shall make this precise in section \ref{physics}, after summarizing the relevant theory underlying equation (\ref{ttt}) in section \ref{results}. The key point is that this theory exists first at the Lie algebra level (step 1), then the representation $\theta$ gives corresponding results for matrices (step 2).
As a preview we just mention here two salient points.  First, it is the asymptotic data at $t=0$ (of a specific solution) which shows that the quantum cohomology of $Gr_k(\C^{n+1})$ is the $k$-th exterior power of the quantum cohomology of $\C P^n$. Second, it is the asymptotic data at $t=\infty$ which shows that the soliton data (or monodromy data) of $Gr_k(\C^{n+1})$ is the $k$-th exterior power of the soliton data (or monodromy data) of $\C P^n$.  The link between the data at zero and the data at infinity is provided by solutions of the tt* equations, as we shall see in the next section.  The existence of these solutions and their explicit asymptotics is nontrivial.

\section{Results on the  tt*-Toda equations}\label{results}

The case $n=1$ of (\ref{ttt}) is the radial sinh-Gordon equation, a special case of the Third Painlev\'e equation.  In fundamental work, motivated by the Ising model, McCoy-Tracy-Wu studied this case in detail, and their results provided an essential test case for the conjectures of Cecotti and Vafa (see \cite{CeVa91},\cite{CeVa92}).  The \ll separatrix solution\rr of the tt* equations in that case was observed to correspond to the quantum cohomology of $\C P^1$. In \cite{Za96}, Zaslow described the prospects for further examples as follows: 

{\em Other spaces, such as the higher projective spaces and Grassmannians, are too unwieldy for a direct analysis. Too little is known about the solutions to the tt* equations, which for 
$\C P^n$ correspond to affine Toda equations. Perhaps the proposed relation between math and physics is best borne out by a rigorous analysis of these equations and their asymptotic properties.}

For the tt*-Toda equations we present results of this kind in this section.
We suppress some technical details, which may be found in the references given.

\begin{theorem}\cite{GuLi14}\cite{GIL1}\cite{GIL2}\cite{MoXX}\label{thm1}
There is a bijection between solutions $(w_0,\dots,w_n)$ of the $tt^*$-Toda equations on $\C^\ast$ and points of the region
\[
\{ (m_0,\dots,m_n)\in\R^{n+1}\mid
m_{i+1}-m_i\leq 1, m_i+m_{n-i}=0\}
\]
in  $\R^{n+1}$.
The correspondence is given by $w_i(t)\sim -m_i\log |t|$, as $|t|\to 0$.
\end{theorem}
The proof in \cite{GuLi14},\cite{GIL1},\cite{GIL2} is by p.d.e.\ methods.  The signs on the right hand side of (\ref{ttt}) facilitate the use of the Maximum Principle, and ensure the uniqueness of the solution with the given asymptotic data.  The method used by Mochizuki in \cite{MoXX} is also based on p.d.e.\ theory, though it is more sophisticated, involving a transformation to a (suitably extended) version of the non-abelian Hodge Correspondence.

To investigate the properties of these solutions we use integrable systems methods and the connection form $\al=\al(t,\bar t)$. 
For this, it will be convenient to introduce a new variable $z$, related to $t$ by
(\ref{ztot}) below.
A key role is played by a {\em holomorphic} (in $z$) connection form $\om=\om(z)$, given by
\begin{equation}\label{momega}
\om=\tfrac1\la\eta\, dz,\quad
\eta=
\bp
 & & & z^{k_0}\\
 z^{k_1} & & & \\
  & \ddots & & \\
   & & z^{k_n} &
   \ep
\end{equation}
where the $k_i$ are real numbers with $k_i\ge -1$, such that $k_i=k_{n-i+1}$ for $1\le i \le n$. If all $k_i\in\Z$ the connection is holomorphic on $\C^\ast$; if some $k_i\notin\Z$ the connection is holomorphic only on the universal covering of $\C^\ast$. Using this holomorphic data we shall construct a connection form $\al$, and hence a family of local solutions of (\ref{ttt}), which will include the global solutions of Theorem \ref{thm1}. 

The relation between $\om$ and $\al$ is explained in sections 2 and 7 of \cite{GIL3}.  We summarize this briefly. Given $\om$, there exist a loop $\ga:S^1\to \SL_{n+1}\C$ and a gauge transformation $G$ such that
\begin{equation}\label{gaugeG}
\al = [ (\ga L)_\R G ]^{-1} d[(\ga L)_\R  G ],
\end{equation}
where $L$ is a normalized solution of the o.d.e.\ 
$
L^{-1}L_z=\tfrac1\la\eta
$
and $(\ga L)_\R$ is the first factor in the loop group Iwasawa factorization
\[
\ga L = (\ga L)_\R (\ga L)_+.
\]
This factorization exists in a neighbourhood of $z=0$. 
Iwasawa factorization is with respect to a real form isomorphic to the loop group of $SL_{n+1}\R$. 
All functions here are in general multi-valued on $\C^\ast$.

The normalized $L$ is of the form
\[
L(z,\la)=e^{\frac1\la \cN \log z}\left( I + \sum_{i\ge1} \la^{-i} S_i(z)\right)
\]
where $\cN$ is nilpotent and $S_i(0)=0$ for all $i$. It takes values in the
loop group of $\SL_{n+1}\C$. 
The relation between the variables $t$ and $z$ is: 
\begin{equation}\label{ztot}
\tfrac N{n+1}t=z^{\frac N{n+1}},
\end{equation}
where $N=n+1+\sum_{i=0}^n k_i$.  We exclude the trivial case $N=0$ (i.e.\ all $k_i=-1$). 
The $\SL_{n+1}\C$-valued gauge transformation $G$ is
$G(z,\zbar)=(\vert t\vert/t)^m=(\vert z\vert/z)^{ \frac N{n+1} m}$,
where the matrix $m=\diag(m_0,\dots,m_n)$ is given by $1-m_i+m_{i-1}=\tfrac{n+1}N(k_i+1)$.
Thus $G$ is determined by $k_0,\dots,k_n$. There is some freedom in the choice of $\ga$; this freedom is equivalent to replacing $z^{k_i}$ by $c_i z^{k_i}$ and varying the $c_i \, (>0)$. In \cite{GIL3} a fixed $\ga$ was used and the $c_i$ were varied. Here we fix all $c_i=1$ and allow $\ga$ to vary.

Thus, for each $(k_0,\dots,k_n)$ --- i.e.\ for the holomorphic data $\om$ --- we obtain a solution $w$ of the p.d.e.\ (\ref{ttt}).  By construction it satisfies $w\sim -m\log \vert t\vert$ as $t\to 0$.
However, this does not qualify as a solution of the tt*-Toda equations until it has been proved that the solution is defined on all of $\C^\ast$. In fact the above construction gives a solution on some open set $0<\vert t\vert< \epsilon$, where $\epsilon$ depends on $(k_0,\dots,k_n)$ and $\ga$, and all such solutions \ll near $0$\rr arise this way. Theorem \ref{thm1} ensures that, for each $(k_0,\dots,k_n)$,
 there is precisely one $\ga$ for which the solution is global, i.e.\ for which $\epsilon=\infty$:

\begin{theorem}\cite{GIL3}\cite{GILX}\label{thm2} Fix $N=n+1+\sum_{i=0}^n k_i>0$. 
For each $(k_0,\dots,k_n)$ (with $k_i=k_{n-i+1}$ and $k_i\ge-1$ for all $i$) there exists a unique solution $(w_0,\dots,w_n)$ of (\ref{ttt}) on $\C^\ast$, and all solutions on $\C^\ast$ arise this way. The corresponding point $(m_0,\dots,m_n)$ in Theorem \ref{thm1} is determined uniquely by the equations
$1-m_i+m_{i-1}=\tfrac{n+1}N(k_i+1)$.
\end{theorem}

The fact that the solutions obtained by p.d.e.\ methods in Theorem \ref{thm1} can be constructed from the holomorphic data $\om$ is close to the Hitchin-Kobayashi point of view, if $\om$ is regarded as a Higgs field, and we shall give a precise relation in section \ref{morett*}.

\begin{theorem}\cite{GIL2},\cite{GILX}\label{thm3}
The solution $(w_0,\dots,w_n)$ corresponding to 
$(m_0,\dots,m_n)$ in Theorem \ref{thm1} has (and is uniquely defined by) 
the following asymptotic behaviour as $|t|\to \infty$:  
\begin{equation}\label{FTasympinfinity}
-\tfrac 4{n+1}
\sum_{p=0}^{[\frac12(n-1)]} w_p \sin \tfrac{(2p+1)k\pi}{n+1}  
\sim
 s_k\, F(L_k x),
\quad
k=1,\dots,[\tfrac12(n+1)]
\end{equation}
where 
\[
F(x)=\tfrac12(\pi x)^{-\frac12}e^{-2x}, 
\quad
L_k=2\sin \tfrac k{n+1}\pi,
\]
and $s_k$ is the $k$-th elementary symmetric function of the $n+1$ quantities
\[
e^{ (2m_0+n)\frac{\pi\i}{n+1}  },
e^{ (2m_1+n-2)\frac{\pi\i}{n+1}  },
\dots,
e^{ (2m_n-n)\frac{\pi\i}{n+1}  }.
\]
The notation $[\tfrac12(n+1)]$ means $\tfrac12(n+1)$ if $n$ is odd, and
$\tfrac12n$ if $n$ is even.
\end{theorem}

\no Note that, as $w_i+w_{n-i}=0$, 
the independent functions are $w_0,\dots,w_{[\frac12(n-1)]}$, whether $n$ is odd or even.

\begin{remark} 
The linear system in formula (\ref{FTasympinfinity}) can be \ll solved\rr for $w_p$. When $n$ is odd, for example, the result is
\begin{equation}\label{asympinfinity}
w_p
\sim 
\sum_{k=1}^{\frac12(n-1)} s_k \,  F(L_k x)  \sin \tfrac{(2p+1)k\pi}{n+1}  
\ +\ 
\tfrac12 (-1)^p s_{\frac12(n+1)} F(L_{\frac12(n+1)} x).
\end{equation}
Note that the \ll leading term\rr here depends on which (if any) of the $s_k$ are zero.
\end{remark}

The proof of Theorem \ref{thm3} is superficially analogous to that of Theorem \ref{thm2}: first, local solutions near $t=\infty$ are constructed by solving a Riemann-Hilbert problem, then Theorem \ref{thm1} is invoked to guarantee a global solution. However, there is much more involved, as the setting of the Riemann-Hilbert problem in the case of Theorem \ref{thm3}
requires a link between the asymptotic data at $0$ and $\infty$, and this depends on the isomonodromy formulation of (\ref{ttt}). That is, it depends on the fact that (\ref{ttt}) is the condition for the monodromy data (monodromy matrices, formal monodromy matrices, Stokes matrices, and connection matrices) of the {\em meromorphic} (in $\la$) connection form
\begin{equation}\label{malphahat}
\hat\al=
\left[
- \tfrac{t}{\la^2}
\ W^T
- \tfrac1\la  x w_{x} + \tbar \,W
\right]
d\la,\quad
x=\vert t\vert
\end{equation}
to be independent of the parameter $t$.  
As was pointed out by Dubrovin in \cite{Du93}, 
the existence of such an isomonodromy formulation for the tt* equations follows from a homogeneity condition 
(the radial condition, in the tt*-Toda situation). 

The connection form $\hat\al$ has (in the variable $\la$) poles of order $2$ at $\la=0$ and $\la=\infty$.
The parameter $s_i$ arises in the computation of the Stokes matrices at these (irregular) poles. We explain this briefly next, with more details in section \ref{physics}.

Stokes matrices (as formulated in classical o.d.e.\ theory) generally depend on various non-canonical choices. However, in the case of $\hat\al$, they are equivalent to the data 
\[
\cS=(s_1,\dots,s_n),\quad s_i=s_{n-i+1}
\]
which can be defined canonically as the coefficients of the characteristic polynomial of a certain  matrix $M\in SL_{n+1}\C$.  This $M$ (we shall give a precise definition below) represents (plus or minus) the \ll$\frac1{n+1}$-th monodromy\rr of the connection form $\hat\al$. It can be shown (see sections 3 and 4 of \cite{GH1}) that all Stokes matrices and the monodromy matrix can be recovered from $M$.  In the context of Theorems \ref{thm2} and \ref{thm3}, the formal monodromy is trivial, and the connection matrix will take a fixed canonical value, 
so only the Stokes and monodromy matrices can vary.  Because of this we shall refer to $M$ simply as the {\em monodromy data} from now on. 

\begin{example}\label{cpn} (Complex projective space)
The point $(m_0,\dots,m_n)=(-\tfrac n2,-\tfrac n2+1,\dots,\tfrac n2)$,
i.e.\ $k_0=0$ and $k_1=\cdots =k_n=-1$,  corresponds to the connection form
\[
\om=\tfrac1\la
\bp
 & & & z\\
1 & & & \\
  & \ddots & & \\
   & & 1 &
   \ep
\, \tfrac{dz}{z}
\]
which represents the Dubrovin connection (or quantum D-module) associated to the quantum cohomology of complex projective space $\C P^n$. We refer to
\cite{CoKa99} or \cite{Gu08} for more information on this well known fact. By the last assertion of Theorem \ref{thm3}, for the corresponding solution of (\ref{ttt}), we have
$
s_i=\tbinom{n+1}{i}.  
$
\qed
\end{example}  

To go further (in particular, to explain the role of the Satake Correspondence), it will be useful to formulate the holomorphic data $\om$ and the monodromy data $M$ Lie-theoretically. We refer to sections 5 and 6 of \cite{GH2} for the details of what follows.  

Let $G$ be the (complex, simple) simply-connected Lie group with Lie algebra $\g$. 
First, we take
\begin{equation}\label{gomega}
\om=
\tfrac1\lambda 
\sum_{i=0}^l z^{k_i} e_{-\be_i}
dz
\end{equation}
as the Lie-theoretic version of (\ref{momega}).
The Lie-theoretic
version of $\hat\al$ in (\ref{malphahat}) is
\begin{equation}\label{galphahat}
\hat\al=
\left[
-\tfrac t{\lambda^2} \Ad(e^w) E_-
- \tfrac1\la x w_{ x } 
+ \bar t \Ad(e^{-w}) E_+
\right]
d\la
\end{equation}
where $x=\vert t\vert$. The role of the $i$-th
elementary symmetric function is played by the character $\chi_i$ of
the $i$-th basic irreducible representation of $G$.
We write $\chi=(\chi_1,\dots,\chi_l):G\to\C^l$.
We shall also need the Coxeter number $h=\sum_{i=0}^l q_i$ of $\g$; here we recall that 
the $q_i$ are defined by $\psi=\sum_{i=1}^l q_i\be_i$, and we put $q_0=1$.  

Then (Theorem 6.8 of \cite{GH2}), if $w$ is a (local) solution constructed from the holomorphic data
$\om$, the Stokes data of $\hat\al$ is given by
\begin{equation}\label{gS}
\cS=(s_1,\dots,s_l)=\chi(e^{2\pi\i(m+x_0)/h}),
\end{equation}
where $m\in\h$ is defined by $\be_i(m)=-1+\tfrac h N(k_i+1)$, $1\le i\le l$, and
$N=h+\sum_{i=0}^l q_i k_i$ (as in the case $\g=\sl_{n+1}\C$ we have
$w(\vert t\vert)\sim -m\log |t|$, when $|t|\to 0$). 
The element $x_0\in\h$ is defined by
$x_0=\sum_{i=1}^l \eps_i$, where
$\be_i(\eps_j)=\de_{ij}$.    
The Lie-theoretic version of the monodromy data $M$  is then
\begin{equation}\label{gM}
M=C^\Ga(\cS),
\end{equation}
where $C^\Ga$ is a suitable \ll Steinberg cross-section\rr of $\chi$.

For the definition of $C^\Ga$, see \cite{St74}, Theorem 4, page 120. Explicitly,
for any choice of Cartan subalgebra $\c$ and (ordered set of) simple roots $\Ga=(\phi_1,\dots,\phi_l)$ with respect to $\c$, $C^\Ga$ is defined by
\[
C^\Ga(t_1,\dots,t_l)=\overset{\scriptstyle l}{\underset{i=1}\Pi}
 \exp(t_i e_{\phi_i}) \, n_i,
\]
where $n_i$ is a representative (in $G$) of the 
Weyl group element given by reflection in $\Ker \phi_i$. It is a
cross-section of the map $\chi:G^{\text{reg}}\to\C^l$, where $G^{\text{reg}}$
is the set of regular elements of $G$.

Thus we have the following purely Lie-theoretic statement:

\begin{theorem}\cite{GH2}\label{thm4} 
From the holomorphic data $\om=
\tfrac1\lambda 
\sum_{i=0}^l z^{k_i} e_{-\be_i}
dz$ 
(with $k_i\ge-1$ for all $i$)
we obtain the explicit
monodromy data 
$M=C^\Ga(\cS)=C^\Ga(\chi(e^{2\pi\i(m+x_0)/h}))$.
\end{theorem}

Now let $\th:\g\to\End(V)$ be a faithful representation of $\g$. We denote by
$\Th:G\to\Aut(V)$ the corresponding  representation of $G$. We define
\[
\om_\th=\th(\om)=
\tfrac1\lambda
\sum_{i=0}^l z^{k_i} \th(e_{-\be_i})
dz,
\quad
M_\Th=\Th(M).
\]
We remark that there is also a Lie-theoretic formula for individual Stokes matrices --- to be given later as (\ref{gQ}) in 
section \ref{particles} --- and $\Th$ may be applied to that, in the same way.

\begin{corollary}\label{cor5} 
From the holomorphic data $\om_\th$ we obtain the explicit monodromy data $M_\Th$.
\end{corollary}

Although this corollary is merely a \ll matrix representation\rr of Theorem \ref{thm4}, 
we consider it important as it reveals the functorial nature of the relation between $\om_\theta$ and $M_\Theta$.  Formula (\ref{gM}), the formula for the monodromy data at the Lie group level, is the key ingredient. 

Now we return to the case $\g=\sl_{n+1}\C$, where
we have Theorems \ref{thm1}, \ref{thm2}, \ref{thm3} on the global solutions of the tt*-Toda equations (\ref{ttt}). Combining the results of this section gives: 

\begin{theorem} 
Let $\g=\sl_{n+1}\C$, $G=\SL_{n+1}\C$.
Consider the global solutions of the tt*-Toda equations (\ref{ttt}) 
as in Theorems \ref{thm1}, \ref{thm2}, \ref{thm3}.
Each such solution $w$ gives rise to a correspondence between 
the holomorphic data  $\om_\th$ (constructed from the asymptotic data of $w$ at $t=0$) and the 
monodromy data $M_\Th$ (constructed from the asymptotic data of $w$ at $t=\infty$).
\end{theorem} 

The deceptively simple condition $k_i\ge-1$, which is equivalent to $m_i-m_{i-1}\ge1$, looks more complicated when written in terms of $\cS$ (an example can be seen in section 5 of \cite{GIL2}).  It is analogous to the stability conditions which are needed in the Hitchin-Kobayashi Correspondence. 
Although the formula for $M$ in Theorem \ref{thm4} is meaningful without the restriction $k_i\ge-1$,  it would not (in the absence of this restriction) represent the 
correct data for the connection $\hat\al$; it would not bear any relation to the tt*-Toda p.d.e.\

\begin{example} (Complex Grassmannians)
Consider again the solution of (\ref{ttt}) corresponding to the point
$m=(-\tfrac n2,-\tfrac n2+1,\dots,\tfrac n2)=-x_0$.
For $\theta=\la_{n+1}$, $\om_\theta$ represents
the quantum cohomology ring of $\C P^n$ (as in Example \ref{cpn}),
and $M_\Theta$ represents the Stokes/monodromy matrices of the Dubrovin connection for  $\C P^n$. 
 For $\theta=\wedge^k\la_{n+1}$, $\om_\theta$ represents
the quantum cohomology ring of $Gr_k(\C^{n+1})$ (see section 7 of \cite{CDGXX},
and also \cite{GoMaXX}, \cite{KaXX}), and  $M_\Theta$ represents the Stokes/monodromy matrices of the Dubrovin connection for $Gr_k(\C^{n+1})$.
The canonical Stokes data
$s_i=\tbinom{n+1}{i}$ is the same, for any $k$. 

This shows that the the tt* metrics associated to the quantum cohomology of 
 $\C P^n$ and $Gr_k(\C^{n+1})$ have the same analytic origin: the global solution of (\ref{ttt}) corresponding to the point $m=-x_0$. The \ll initial data\rr (asymptotic data at $t=0$) of the solution is the holomorphic connection form $\om = \frac1\la(z e_{-\be_0} + \sum_{i=1}^n e_{-\be_i})\frac{dz}{z}$.  This connection form has been discovered empirically by other authors --- it is the \ll Frenkel-Gross connection\rr of \cite{FG09}, and was an important ingredient of \cite{GoMaXX}, where the link with Kostant's work \cite{Ko59} was attributed to
 \cite{GoPe12}; the latter was apparently motivated by work of Lerche-Vafa-Warner which we shall discuss in the next section. We emphasize that this particular connection form  {\em originates} naturally from the tt*-Toda equations. 
 
The (uniqueness of the) global solution of (\ref{ttt}) corresponding to $m=-x_0$ also has practical consequences: it leads to a canonical choice of fundamental solution of the quantum differential equation (the equation for flat sections of $d+\om$), and hence a canonical choice of \ll central connection matrix\rr (cf.\ \cite{GIL3}).
\qed
\end{example}

\section{Applications to field theories}\label{physics}

We return now to the field theory point of view.  Our aim is to explain how the results of section \ref{results} provide a mathematical foundation for the article \cite{Bo95a} of Bourdeau and for related results in the physics literature.  For this we need to describe some of the physics aspects more precisely. 

\subsection{The physical data}\label{morephysics}

The physical data at $t=0$ is the chiral ring structure on the vector space $V$ of states.   A state $\vert\al\rangle$ corresponds to a chiral operator $\phi$, where $\vert\al\rangle=\phi\,\vert0\rangle$. Composition of operators gives a product operation on $V$. The matrix of multiplication by the $i$-th basis element of $V$ is denoted by $C_i$.
The structure constants (with respect to the basis) are the $3$-point correlation functions; there is also a $\C$-linear inner product on $V$, the topological metric, and its entries are the $2$-point correlation functions. In the context of quantum cohomology, these are the Gromov-Witten invariants and intersection form, respectively. 

Regarded as in the Neveu-Schwarz sector, each state has a chiral charge and there is a distinguished state of charge $0$ (essentially the cohomology grading, in the case of quantum cohomology). Regarded as in Ramond sector, the charges are normalized to be symmetric about $0$. 

The chiral operators are themselves holomorphic functions of $V$;  more precisely, they are endomorphisms of the tangent bundle of $V$. 
They satisfy the associativity equations (WDVV equations). The data described so far has been abstracted in the mathematical theory of Frobenius manifolds.  Physically, this holomorphic data is regarded as \ll topological\rrr.  

The tt* equations (see section \ref{morett*}) are the equations for pluriharmonic metrics $g$ associated to the holomorphic data. Such metrics combine the holomorphic data with complex conjugate antiholomorphic data; this is the \ll fusion\rr of topological data with antitopological data.  It has also been abstracted, in the mathematical theory of harmonic bundles/Higgs bundles.

In our tt*-Toda example, $V=\C^{n+1}$.  
The topological metric is represented by the anti-diagonal matrix $\De=(\de_{i,n-i})_{0\le i\le n}$. The Neveu-Schwarz  charges are (related to) our $k_0+1,\dots,k_n+1$ and the Ramond charges are our $m_0,\dots,m_n$. Either of these parametrize the solutions of the tt*-Toda equations (though only for certain special parameter values, for example $m=-x_0$,  is there a quantum cohomology interpretation).

In our situation, and in the examples of \cite{CeVa91},\cite{CeVa92},
the product operation depends only on a (complex) one-dimensional subspace of $V$ (this means that only \ll small quantum cohomology\rr is being treated). 
Instead of the tangent bundle of $V$, we are working with its restriction to the one-dimensional subspace, and we just have one chiral matrix $C$.
The WDVV part of the tt* equations becomes trivial; we have harmonic maps rather than general pluriharmonic maps. 
Some vestiges of Frobenius manifold theory remain, such as the homogeneity (Euler) condition and the topological metric, but the harmonic/Higgs theory is now the main focus.  

The physical data at $t=\infty$ consists of Bogomolnyi solitons interpolating between vacua; this arises from path-integral considerations (see pages 345-6 of \cite{Za96}), and in our situation is detected in the leading term of the asymptotics of a solution $w$ as $t\to\infty$.  
The multiplicity of the $i$-th soliton is the $s_i$ of Theorem \ref{thm3}.  
For this reason  $s_i$ must be a nonnegative integer
if $w$ corresponds to a physically realistic model.
Which vacua are relevant for the $i$-th soliton will be explained later.

\subsection{The topological-antitopological fusion equations}\label{morett*}

The tt* equations, as formulated by Cecotti and Vafa 
in section 3 of \cite{CeVa91}, are
\begin{equation}\label{ceva}
\bar\b(g\b g^{-1}) -  [C,g C^\dagger g^{-1}]=0
\quad
(\b=\tfrac{\b}{\b t}, \bar\b=\tfrac{\b}{\b \bar t})
\end{equation}
where $C$ is the (holomorphic) chiral matrix, $C^\dagger=\bar C^T$, and $g^{-1}$ is the Hermitian matrix representing the tt* metric, i.e.\ the metric is $(x,y)\mapsto \bar x^T g^{-1} y$.  
Thus $g C^\dagger g^{-1}$
is the Hermitian adjoint of $C$ with respect to $g^{-1}$. 
This
is equation (3.9) in \cite{CeVa91}; the companion equation (3.10) there is vacuous as we assume only one chiral matrix.   All matrices here are defined with respect to a holomorphic local frame $h_0,\dots,h_n$, or \ll holomorphic gauge\rrr.

Section 4.2 of \cite{CeVa92} gives a zero-curvature formulation 
\begin{equation}\label{zcchol}
[\b + g\b g^{-1} -xC,\bar\b - \tfrac1x g\,C^\dagger g^{-1}]=0
\end{equation}
of (\ref{ceva}), due to Dubrovin, where $x$ is a parameter. Expanding in powers of $x$ and setting the coefficients equal to zero, we see that it is equivalent to (\ref{ceva}) together with $\bar\b C=0$. 
The link with harmonic/Higgs theory is now visible: (\ref{zcchol}) has the form of the Higgs bundle equations
\[
[\b + h^{-1}\b h+x\Theta^\pr,\bar\b + \tfrac1x \Theta^\prr]=0
\]
where $h$ is a Hermitian metric and $\Theta = \Theta^\pr dt + \Theta^\prr d\bar t$ is a Higgs field. 
For this reason we have called $g^{-1}$ the tt* metric, rather than $g$. 

A more symmetrical form is
\begin{equation}\label{zcc}
[\b - \b e \, e^{-1} - x e^{-1} C e, \bar\b + e^{-1} \bar\b e - \tfrac1x e C^\dagger e^{-1}]=0
\end{equation}
where $e=\sqrt g$. 
The gauge transformation
\[
\b+X\mapsto \b+e\,\b e^{-1} + eXe^{-1},\quad
\bar\b+Y\mapsto \bar\b+e\,\bar\b e^{-1} + eYe^{-1}
\]
by $e$
converts (\ref{zcc}) to (\ref{zcchol}).
The symmetrical zero-curvature condition (\ref{zcc}) is equivalent to the tt* equation (\ref{ceva}) together with
the condition
\[
\bar \b(e^{-1}Ce)+[e^{-1}\bar \b e, e^{-1}Ce]=0,
\]
which says that $e^{-1}Ce$ is holomorphic (as an endomorphism) with respect to the $\bar\b$-operator
$\bar\b + e^{-1}\bar\b e$.  The gauge transformation $e$ converts $\bar\b + e^{-1}\bar\b e$
to $\bar\b$, and $e^{-1}Ce$ to $C$.  

Let us denote the new \ll symmetrical\rr frame by $e_0,\dots,e_n$,
so that $e_i=e h_i$.  In this frame the chiral matrix is $e^{-1}Ce$ and the tt* metric is
$\bar e^T g^{-1} e = I$, i.e. just the standard Hermitian metric (thus, the frame is orthonormal).

Compatibility with the topological metric $(x,y)\mapsto  x^T \De y$
imposes further conditions on (\ref{zcchol}), (\ref{zcc}).
In the holomorphic frame these conditions are $\De C^T\De=C$ and $\De g^{-T} \De = g$, where $g^{-T}$ denotes the transpose of the matrix $g^{-1}$.
The latter condition means that $g$ (and also $e$) belong to the $\De$-modified complex orthogonal group
\[
\SO^\De_{n+1}\C = \{ A\in \SL_{n+1}\C \st \De A^{-T} \De = A \}.
\]
Then our flat connection is an $\sl^\De_{n+1}\R$-valued connection, where
\[
\SL^\De_{n+1}\R = \{ A\in \SL_{n+1}\C \st \De \bar A \De = A \}.
\]
Let $F$ be a (locally defined) fundamental solution matrix of the corresponding linear system
\[
\b F =F( - \b e \, e^{-1} - x e^{-1} C e),\ \ 
\bar\b F=F( e^{-1} \bar\b e - \tfrac1x e C^\dagger e^{-1}).
\]
Then (see \cite{Gu08}, section 7.4) the zero-curvature condition (\ref{zcc}) is exactly the
condition that the  map $F\vert_{x=1}$ 
induces a harmonic map to $\SL^\De_{n+1}\C / \SO^\De_{n+1}$.

We have already seen that a solution of the tt* equations defines a Higgs bundle.  It is well known that this is equivalent to being a harmonic bundle, i.e.\ (locally) a 
harmonic map to $\SL_{n+1}\C/\SU_{n+1}$.  Compatibility with the topological metric means that the image of the map lies in the smaller space $\SL^\De_{n+1}\R / \SO^\De_{n+1}$, which is (see
\cite{GuLi14}, section 2) 
isomorphic to $\SL_{n+1}\R / \SO_{n+1}$, the space of all inner products on $\R^{n+1}$. For
this reason the topological metric condition is sometimes referred to as a \ll reality condition\rr (although the symmetric space $\SL_{n+1}\C/\SU_{n+1}$ is
already real).

The relation with Hitchin's work (\cite{Hi92} and references therein) was mentioned incidentally\footnote{(see the footnote on page 371 of in \cite{CeVa91})
} in the original work of Cecotti and Vafa. In more recent work on tt* geometry (notably \cite{GMN13},\cite{CGV14}),
the Higgs point of view has been emphasized.

Our tt*-Toda equation is the following special case of (\ref{zcc})
\begin{equation}\label{zccttt}
[\b + w_t +\tfrac1\la e^{w}E_- e^{-w},  \bar\b  -w_{\bar t} + \la e^{-w} E_+ e^{w}]=0,
\end{equation}
in which we have
\[
C=E_-, \ C^\dagger=E_+,\  g=e^{-2w},\ e=e^{-w},\ x=-1/\la.
\]
The corresponding special case of (\ref{zcchol}) is
\begin{equation}\label{zccholttt}
[\b + 2w_t +\tfrac1\la E_-,  \bar\b + \la g E_+ g^{-1}]=0.
\end{equation}
The holomorphic chiral matrix $C$
is constant here, but we claim that it corresponds to our original \ll Higgs field\rr $\om$ from (\ref{momega})
through a holomorphic gauge transformation and
change of variable from $t$ to $z$.  

To see this, we apply to (\ref{zccholttt})
a holomorphic
gauge transformation of the form $k=at^m$ where $a=\diag(a_0,\dots,a_n)$ is to be chosen
later.  
The gauged connection form is
\begin{equation}\label{gauged}
( k^{-1}k_t + 2w_t +\tfrac1\la k^{-1} E_- \, k)dt+
( \la k^{-1} g E_+ \, g^{-1} k )d\bar t.
\end{equation}
The chiral matrix part $\tfrac1\la k^{-1} E_- k\, dt$ of this is
\[
\tfrac1\la
\bp
 & & & \frac{a_n}{a_0}\,t^{m_n-m_0+1}\\
\frac{a_0}{a_1}\,  t^{m_0-m_1+1} & & & \\
  & \ddots & & \\
   & & \frac{a_{n-1}}{a_n}\, t^{m_{n-1}-m_n+1} &
\ep
\tfrac{dt}{t}.  
\]
From (\ref{ztot}) we have 
$\tfrac N{n+1}t=z^{\frac N{n+1}}$, hence $\frac{dt}t=\frac{N}{n+1}\frac{dz}{z}$. From
Theorem \ref{thm2} we have
$m_{i-1}-m_i+1=\tfrac{n+1}N(k_i+1)$. Converting to $z$ gives exactly
$\om$ if we choose the $a_i$ so that
$a_i a_{n-i}=1$ and 
$ a_{i-1}/a_i = \frac{n+1}{N} ( \frac N{n+1} ) ^{  \frac{n+1}{N} (k_i+1) }$.
Thus our holomorphic data $\om$ is indeed equivalent to the chiral matrix and to the Higgs field.

\begin{remark} (Gaiotto's Conjecture)  The remaining terms of (\ref{gauged})
can be \ll removed\rr by taking a suitable limit, leaving just the chiral (Higgs) part. This gives a direct computational relation between $\al$ and $\om$.
To verify
this, we put $s=t/\la$ and let $t,\la\to 0$ without changing $s$. The $dt$ term
is $(2tw_t+m)\frac{ds}s$ and the nonzero entries in the matrix of the $d\bar t$ term
are of the form 
$\la \bar t e^{-2w_{i-1}+2w_i} t^{-m_{i-1}+m_n+1} \frac{d\bar s}{\bar s}$.
Now we use the fact that $w_i\sim -m_i\log\vert t\vert$ as $t\to 0$, and
$m_{i-1}-m_i+1=\tfrac{n+1}N(k_i+1)\ge 0$. The limit of each term is zero.
Thus, the (suitably scaled) limit of (\ref{gauged}) is just $\om$.
This fact, a special case of a conjecture by Gaiotto in \cite{GaXX}, is thus a
consequence of the Iwasawa factorization construction of local solutions of the tt*-Toda
equations.
\end{remark}

\subsection{Particles and the Coxeter Plane}\label{particles}

Having explained how our tt*-Toda equations fit into the physics framework, we can now give some further applications of the results of section \ref{results}.  Although we are mainly concerned with the $A_n$ case in this article, it will be conceptually clearer to continue using the Lie-theoretic notation introduced in earlier sections. 

First we need the relation between particles in affine Toda field theory and roots, which was proposed in
\cite{BCDS90},\cite{Do91},\cite{Do92},\cite{Fr91}. 
Recall that we have chosen a Cartan subalgebra $\h$ and root system $\De$, and simple roots $\be_1,\dots,\be_l$. Let us choose also a Coxeter element $\ga$ (the product of the reflections in the simple roots).   This is an element of order $h$ in the Weyl group of $\g$, and it acts naturally on $\De$. It is well known (see \cite{Ko59}) that $\De$ decomposes into $l$ orbits, each consisting of $h$ roots, and also that the eigenvalues of $\ga$ on $\h^\ast$ are $e^{2\pi\i n_1/h},\dots,e^{2\pi\i n_l/h}$ where
$1= n_1\le \cdots \le n_l=h-1$ are the exponents of $\g$. 

The first eigenvalue $e^{2\pi\i/h}$ always has multiplicity $1$, and the corresponding eigenspace determines a two-dimensional real plane in $\h^\ast_\sharp$, called the Coxeter Plane.  
(The complexification of this real plane is the sum of the eigenspaces of $e^{\pm2\pi\i/h}$.) 
Projecting the $hl$ roots (regarded as elements of $\h^\ast_\sharp$) orthogonally with respect to the Killing form, we obtain $hl$ points in the Coxeter Plane. It is known  (see \cite{Ko59},\cite{Ko10}) that $\ga$ acts on the Coxeter Plane (and hence on this set of $hl$ points) by rotation through the angle $2\pi/h$.

\begin{figure}[h]
\begin{center}
\includegraphics[scale=0.3, trim= 0 200 0 200]{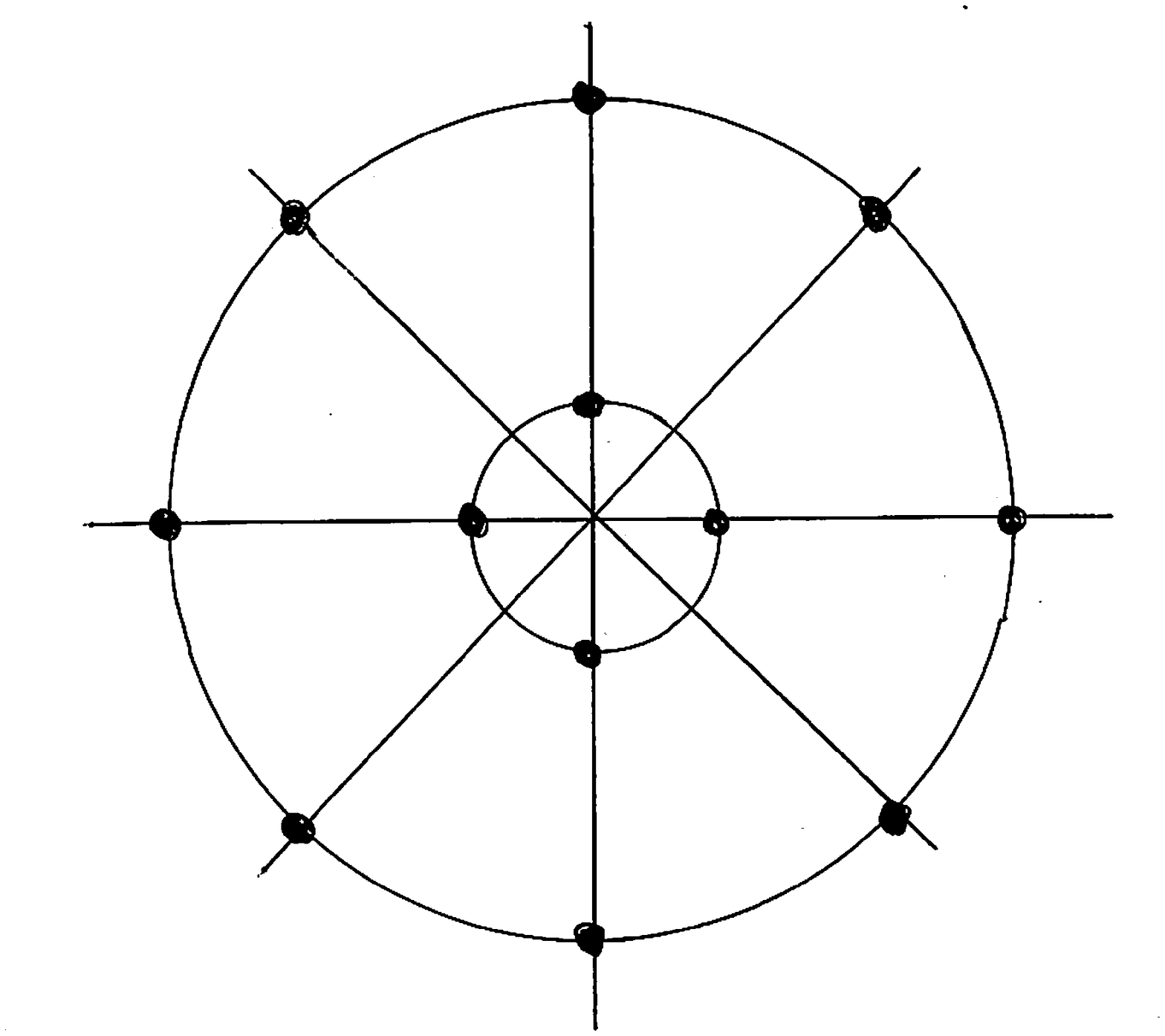}
\end{center}
\caption{Example of a Coxeter Plane.}\label{cox}
\end{figure}

Drawing the rays (\ll spokes\rrr) from the origin to each of these points, and the concentric circles (\ll wheels\rrr) passing through them, we obtain a picture of the type shown in Figure \ref{cox}. 
Each wheel contains an orbit of the action of $\ga$ on $\De$ (in 
some cases, more than one; such orbits may or may not be coincident).

The model proposed in
\cite{BCDS90},\cite{Do91},\cite{Do92},\cite{Fr91}
{\em associates a particle to each orbit (or wheel), the mass of the particle being (up to an overall scale factor) the radius of the wheel.} Thus there are $l$ particles (though some may coincide). Let us denote by $[\be]$ the particle given by the orbit of (the projection of) $\be\in\De$, and by $L_{[\be]}$ its mass. 

Similarly, the $i$-th eigenvalue of $\ga$ gives rise (if it has multiplicity one) to an $i$-th Coxeter Plane, again with spokes and wheels.  In this case the radius of the wheel is declared to be the $n_i$-spin of $[\be]$ (mass is $n_1$-spin, i.e.\  $1$-spin).  

This structure of particles and spins passes various physical consistency checks. The case of $E_8$ has even been investigated experimentally (\cite{BoGa11},\cite{Ko10}). We emphasize that this structure depends only on $G$ (not a particular matrix representation of $G$) and that the tt*-Toda equations did not (yet) play any role in its description.

Nevertheless, as was shown in \cite{GH2}, the particle/mass structure {\em is} directly related to the Stokes Phenomenon for the meromorphic connection form $\hat\al$ of (\ref{galphahat}). Namely, {\em the $2l$ rays in the Coxeter Plane may be identified with the $2l$ singular directions in the $\la$-plane for the pole of order $2$ of the meromorphic connection form $\hat\al$ at $\la=0$.}  In order to explain this, we shall sketch the proof of Theorem \ref{thm3} in Lie-theoretic terms.

To prepare, we recall that singular directions can be defined for $\hat\al$ (at $\la=0$) as the rays (starting at the origin, in the $\la$-plane) which bisect the intersections of successive maximal Stokes sectors, and thus index \ll Stokes factors\rr (see section 3 of \cite{GH1} and section 4 of \cite{GH2}). There are $2l$ such rays.
The products of any $l$ successive Stokes factors are \ll Stokes matrices\rrr.  Thus all Stokes data can be obtained from these Stokes factors.   

As shown in Appendix B of \cite{GH2}, the relation between the singular directions and $\hat\al$ leads to a simple \ll differential equation theoretic\rr  description of the $hl$ points in the Coxeter Plane. In this description the $2$-plane is just $\C$ (rather than a somewhat implicitly defined subspace of $\h_\sharp^\ast$),  and the points are obtained by evaluating the roots on the coefficient of $\la^{-2}$ in
the dual\footnote{The same points would be obtained using $\hat\al$, but we use the dual $-\hat\al^T$ for consistency with \cite{GH2}.}
of $\hat\al$ (rather than by taking their orthogonal projections to that $2$-plane). 

As might be expected there is something to pay for this simplicity:  the roots must be those associated to the unique Cartan subalgebra $\h^\pr$ containing that leading coefficient, which is (up to a scalar factor) $E_+$, rather than the roots associated to $\h$. In the terminology of \cite{Ko59}, the Cartan subalgebra $\h^\pr$ is said to be in apposition to the original Cartan subalgebra $\h$.

The Stokes factor corresponding to a singular ray can be described Lie-theoretically, as follows.  Let us denote by $R(\phi)$ the set of roots $\bed\in\Ded$ (with respect to $\hd$) which project to points on the singular direction with argument $\phi$. Then the Stokes factor corresponding to this ray is of the form
\begin{equation}\label{gQ}
Q_\phi=\exp\left(
\sum_{\bed\in R(\phi)}  s_{\bed} e_{\bed}
\right)
\ \in G
\end{equation}
where the $ e_{\bed}$ are suitably normalized root vectors; here we follow Boalch's formulation in \cite{Bo02}.  The symmetries of the
connection form $\hat\al$ show that $s_{\bed_1}=s_{\bed_2}$ if
$\bed_1,\bed_2$ are in the same Coxeter orbit, that $s_{\bed}=s_{-\bed}$, and that
all $s_{\bed}$ are real. 

In particular, $s_{\bed}$ only depends on $[\bed]$ (the particle), so we may write  $s_{\bed}=s_{[\bed]}$.  The $s_{\bed}$ are (up to permutation) the numbers $s_1,\dots,s_l$ introduced earlier in formula (\ref{gS}).  It follows immediately from our  description of the Coxeter Plane that the radius of the wheel containing $[\bed]$ is
\[
\vert \bed(E_+)\vert = L_{[\bed]}= \text{mass of the particle $[\bed]$}.
\]
Here we ignore the scalar factor in the coefficient
as this overall ambiguity is allowed in the definition of mass,
and we ignore the effect of $\Ad e^{-w}$ because it corresponds to a
harmless conjugation of $\hd$ (which does not change lengths).

Now we can sketch the argument which relates the parameters $s_{[\bed]}$ to the asymptotics as $t\to\infty$ of solutions to the tt*-Toda equations.  For the $A_n$ case full details can be found in \cite{GIL2},\cite{GILX}. 

The idea is to use the singular directions --- the rays of the Coxeter Plane --- as the contour for a Riemann-Hilbert problem, whose solutions (if such exist) would give rise to solutions of the tt*-Toda equations.  Let us denote this contour by $\Ga$, and let us define a piecewise constant \ll jump function\rr
$Q_\Ga:\Ga\to \SL_{n+1}\C$ by
\[
Q_\Ga\,\vert_{ {\ }_\text{ray with argument $\phi$}}=
Q_\phi
\]
where the $ s_{\bed}$ in $Q_\phi$
are now arbitrary real numbers.  The Riemann-Hilbert problem is to find a (piecewise) holomorphic function on $\C^\ast-\Ga$, holomorphic on each of the $2l$ sectors,  whose discontinuities (jumps) are given by $Q_\Ga$. Such a function would (almost tautologically) serve as a fundamental solution matrix of $\hat\al$, hence would (in principal) lead to a solution $w$.

In order to formulate a Riemann-Hilbert problem to which standard methods of solution (as in \cite{FIKN06}, for example) apply, we replace $Q_\Ga$ by
\[
G_\Ga=\left( \exp{\tfrac t\la  E_+} \, \exp{ \tbar\la E_-}\right) Q_\Ga
\left( \exp{\tfrac t\la  E_+} \, \exp{ \tbar\la E_-}\right)^{-1}.
\]
For $\la$ on the ray with argument $\phi$, i.e.,  $\la=ke^{\i\phi}$ with $k>0$,
we have
\begin{align*}
G_\Ga\,\vert_{ {\ }_\text{ray with argument $\phi$}}&=\exp\left(
\sum_{\bed\in R(\phi)}  s_{\bed}   \Ad
(  \exp{\tfrac t\la  E_+} \, \exp{ \tbar\la E_-} ) e_{\bed}
\right)
\\
&=\exp\left(
\sum_{\bed\in R(\phi)}  s_{\bed}  \ 
e^{\tfrac t\la  \ad E_+} \, e^{ \tbar\la \ad E_-} \ 
e_{\bed}
\right)
\\
&=\exp\left(
\sum_{\bed\in R(\phi)}  s_{\bed}  
e^{-x(\frac1k + k)  \vert \bed(E_+)\vert}\
e_{\bed}
\right).
\end{align*}
Here we make use of the fact (from \cite{Ko10}) that $E_-=\bar{E}_+$, 
hence $\vert \bed(E_+)\vert=\vert \bed(E_-)\vert$,
where the bar denotes
conjugation with respect to the real form $\h_\sharp^\pr$ of $\h^\pr$.

We see from this that $G_\Ga$ approaches $I$ (exponentially) on each ray $\ze = k e^{\i\phi}$ as $k\to\infty$ or $k\to 0$.  
By Theorem 8.1 of \cite{FIKN06}),  the Riemann-Hilbert problem --- when formulated as a linear integral equation --- is uniquely solvable for $x=\vert t\vert$ sufficiently large (for any $s_{[\bed]}\in\R$). This means that there exists a
(piecewise) holomorphic function $Y=Y(\ze,x)$  whose discontinuities on the $2l$ rays are given precisely by $G_\Ga$.  
Moreover, this $Y$ extends continuously to $\la=0$, where it must be of the form $e^{2w}$ for some solution $w=w(x)$ of the tt*-Toda equations
(and this solution is smooth for $x$ large, and approaches zero as $x\to\infty$).  

From the integral equation we obtain
\[
\displaystyle Y(0,x) \sim I  + \tfrac{1}{2\pi\ii}\int_{\Ga}  
\frac{G(\ze,x) - I}{\ze} \  d\ze
\]
as $x\to\infty$ (here, for notational convenience,
we are regarding $G$ as a matrix group). 
Integrating, we obtain
\[
Y(0,x)\sim I+
\tfrac{1}{2\pi\ii}
\sum_{\bed\in \Ded}  s_{\bed} 
\int_{0}^\infty  
e^{-x(\frac1k + k)  \vert \bed(E_+)\vert}  \tfrac{dk}k\ 
e_{\bed}.
\]
Laplace's method gives
$\frac1{2\pi}\int_0^\infty  e^{-x(\frac1k+k)}  \frac{dk}k
\sim   \tfrac12 (\pi x)^{-\frac12} \, e^{-2x}$ as $x\to\infty$. Replacing $Y(0,x)$ by $e^{-2w}\sim I-2w$,
we obtain
\[
-2w
\sim
-\i
\sum_{\bed\in \Ded}  s_{\bed} 
F(\vert \bed(E_+)\vert \, x)\
e_{\bed},
\]
where $F(x)=\tfrac12 (\pi x)^{-\frac12} \, e^{-2x}$. 
This can be written as 
\begin{equation}\label{gasymptotics}
w
\sim
\sum_{[\bed]}  
s_{[\bed]} 
F(L_{[\bed]}\, x)\,
\left(
\tfrac12{\i}
\textstyle\sum_{\ga^\pr\in [\bed]}
e_{\ga^\pr}
\right).
\end{equation}
In section \ref{bourdeau} we shall make  (\ref{gasymptotics}) more explicit in the $A_n$ case, and
show how it gives (\ref{FTasympinfinity}) and  (\ref{asympinfinity}). 

Formula (\ref{gasymptotics}) demonstrates, at the very least, that the relation between particles/masses and the Stokes data of $\hat\al$  (via the Coxeter Plane and the singular directions) is not merely superficial.  It shows that the predicted particle/mass data $[\bed]$ and $L_{[\bed]}$ appear naturally in the asymptotics of the solution $w$, and it produces the coefficients $s_{[\bed]}$, which can be interpreted as \ll soliton multiplicities\rrr.  It will be exploited further in
the next section, in certain types of model,
to describe the solitons themselves more explicitly.

\begin{remark}  The operator $\ad E_+ \ad E_-$ ($=\ad E_- \ad E_+$), with eigenvalues
$L_{[\bed]}^2=\vert \bed(E_+)\vert^2$ on $\h$, 
was the central ingredient in the article \cite{BrSc16} by Brillon and Schechtman (also motivated by the particle/mass interpretation of the Toda equations).  It arises there because the linearization of 
\[
(2w_{t \bar t}=)\ 
\tfrac12(w_{xx}+\tfrac1x w_x)=[ \Ad(e^w) E_-,\Ad(e^{-w}) E_+ ]
\]
is just
\[
\tfrac12w_{xx} = [E_+, [E_-,w]]
\]
(assuming rapid decrease of $w$ as $x\to\infty$).
As in our situation, the computation of the eigenvalues uses properties of $\h^\pr$ from \cite{Ko59},\cite{Ko10}, so the Coxeter Plane is implicit in \cite{BrSc16} too.  
On the other hand  \cite{BrSc16} emphasizes the Cartan matrix, rather than the Coxeter element, using a well known relation between the two.  From this viewpoint, the masses arise in \cite{BrSc16} (as they did originally in \cite{Fr91}) as the coordinates of the 
Perron-Frobenius eigenvector of the Cartan matrix.
\end{remark}

\subsection{Solitons in polytopic models} 
Now we come to the  polytopic models of  Fendley, Lerche, Mathur, and Warner 
(\cite{FLMW91},\cite{LeWa91}).  Here we shall make use of a global solution $\th(w)$ of the tt*-Toda equations together with a choice of faithful  representation
$\th:\g\to\End(V)$ of $\g$.  

Let $V=\oplus_{i=1}^N \C v_i$, where $v_i\in V$ is a
weight vector, with corresponding weight $\la_i\in\h^\ast$, and
let us assume for simplicity that all weights have multiplicity one (as
will be the case in our examples). 
In this situation, the solitons --- which we have glimpsed rather indirectly so far, just through their attributes of mass and multiplicity, without actually identifying them as particles tunnelling between vacua --- will be made more concrete.

The proposal of \cite{LeWa91} (section 4) is this:  
{\em there is (fundamental) soliton between vacua $v_i$ and $v_j$ precisely when the difference of the corresponding weights $\la_i-\la_j$ is
a\footnote{(i.e.\ one root, as opposed to a sum of roots)
} root.  If 
$\la_i-\la_j=\be\in\De$ then the soliton is (or is dual to) a particle of type $[\be]$, and
its mass is $L_{[\be]}$.}  

The polytope in $\h^\ast$ whose vertices are the weights of $\th$ is called the soliton polytope.  Its projection to the Coxeter Plane provides a visualization of the solitons and the relevant vacua.  Namely, {\em the mass of a soliton between two given vertices is the length of the straight line between the projected vertices}  (and similarly for the $n_i$-spin in 
the $i$-th Coxeter Plane).

Although Lerche and Vafa make these assertions only as a \ll working hypothesis\rrr, they go on to present \ll compelling evidence\rrr, from the point of view of physical consistency.   
The tt*-Toda equations (whose solutions did not play any role in the theory of Lerche and Vafa) contribute further evidence to this proposal.  The tt*-Toda equations also provide more information:  the Stokes parameter $s_{[\be]}$ (which was absent in \cite{FLMW91},\cite{LeWa91})
can be interpreted as the multiplicity of the particle (soliton) $[\be]$ whose mass is $L_{[\be]}$. This will be made precise in the examples below.

\begin{example} (SLOHSS models \cite{FLMW91},\cite{LeWa91}) 
Although we have used any (faithful) representation $\th$ in the above formulation, the polytopic models of \cite{FLMW91},\cite{LeWa91} require $\th$ to be minuscule, i.e.\ that the weights $\la_1,\dots,\la_N$ form a single orbit of the action of the Weyl group. This is a very strong restriction, which implies that the vector space $V$ can be identified with $H^\ast(G/P;\C)$, the cohomology of the projectivized maximal weight orbit $G/P$ of $G$.  The weight spaces are the Neveu-Schwarz vacua, and the chiral charges are (half of) the cohomological degrees.  Striking pictures of projections of soliton polytopes (to the Coxeter Plane, and to higher  Coxeter Planes) can be found in section 5 of  \cite{LeWa91}. 
From the tt* point of view (at least in the case of Theorem \ref{thm1},  this model corresponds to the global solution given by $m=-x_0$. This is essentially the sigma-model of the K\"ahler manifold $G/P$.  Then $\om_\th$ is the matrix of quantum multiplication 
by the K\"ahler class (for a discussion of this, and further references, see \cite{GoMaXX}, \cite{KaXX}).  The Stokes data is given by (\ref{gS}) as 
\[
s_i=\chi_i(1)=\dim\th_i, 1\le i\le l
\]
where $\chi_i$ is the 
$i$-th basic irreducible representation $\th_i$ of $G$.   These dimensions are, of course, nonnegative integers. All minuscule representations of $G$ give rise to this same set of Stokes data.
\qed
\end{example}

The Grassmannian sigma-model --- our main interest in this article --- is a case to which all of the theory described in this section applies.  As promised in the introduction, we return now to this case and to the article  \cite{Bo95a} of Bourdeau.

\subsection{Topological-antitopological fusion and the quantum cohomology of Grassmannians}\label{bourdeau}

Now we take the Lie algebra $\g=\sl_{n+1}\C$ with the standard diagonal Cartan subalgebra $\h$ and simple roots $x_0-x_1$,\dots,$x_{n-1}-x_n$.  We take  $\theta=\wedge^k\la_{n+1}$.  For $k=1,\dots,n$ these give the minuscule representations.

First we have to clarify the particle/mass structure. This is the same for all global solutions.  
As $l=n$,  there are (a priori) $n$ particles, corresponding to the $n$ orbits of 
 the cyclic permutation $(n\ n\!-\!\!1\  \dots \ 1\ 0)$ on
the roots $\De=\{ x_i-x_j \st i\ne j \}$.  We shall see that some of these orbits coincide after projection to the Coxeter Plane, resulting in
$[\frac12(n+1)]$ distinct particles.

To compute the projections (by our simplified method), 
it is convenient to choose an identification $\h^\pr\cong\h$, so that we may compute using the original roots $\De$ rather than the apposition roots $\Ded$.  This may be accomplished using the Vandermonde matrix
$\Om=(\omega^{ij})_{0\le i,j\le n}$, with $\om=e^{{2\pi\i}/{(n+1)}}$. In the present situation we have
\[
E_+=
E_{n,0}+\sum_{i=0}^{n-1} E_{i,i+1}
%=
%\bp
 %& 1& & \\
% & &\ddots & \\
%  &  & & 1\\
%  1 & &  &
%   \ep
   \in\sl_{n+1}\C
\]
and $E_-=E_+^T$.
We have
$
\Om^{-1}  \h^\pr \Om = \h
$
because $\h^\pr$ is spanned by the matrices $E_+,\dots,E_+^n$ and the columns of $\Om$ are the eigenvectors of $E_+$ with eigenvalues $1,\om,\dots,\om^n$ respectively.
Thus
$
\Om^{-1} E_+ \Om = d_{n+1},
\quad
\Om^{-1} E_- \Om = d_{n+1}^{-1}
$
where $d_{n+1}=\diag(1,\om,\dots,\om^n)$. 

For compatibility with \cite{GH2} we shall evaluate apposition roots on $-E_+$, rather than $E_+$.
The projection of the root $\be=x_i-x_j$ to the Coxeter Plane is then given by
$\be(-d_{n+1})=-(\om^i-\om^j)\in\C$,
and its length (the radius of the wheel on which it lies) is just
\[
\vert \be(-d_{n+1}) \vert =
\vert \om^i-\om^j \vert = 2\sin \tfrac {\vert i-j \vert} {n+1}\pi.
\]
We obtain  $[\frac12(n+1)]$ wheels with these radii, each wheel containing $n+1$ equally-spaced points. 
This means that there are $[\frac12(n+1)]$ distinct particles, and the mass of the $k$-th particle is
$L_k=2\sin \tfrac {k} {n+1}\pi$. 
More generally, from $(-d_{n+1})^r$, the
$r$-spin
of the $k$-th particle is $2\sin \tfrac {rk}{n+1}\pi$.  Thus
we have the explicit particle/mass data. 

Next, we discuss data which depends on the particular model (i.e.\ solution of the tt*-Toda equations).  For the global solution $w$
corresponding to $m$, the Stokes data $\cS=(s_1,\dots,s_n)$  is  given in terms of $m$ by the formula of Theorem \ref{thm3}.
The formula for the asymptotics of $w$ as $x\to\infty$, although stated in (\ref{FTasympinfinity}) and
(\ref{asympinfinity}) for the global solutions, actually applies to any solution which is smooth near infinity, as the Riemann-Hilbert argument sketched above shows.   
Formula (\ref{asympinfinity}) is thus a special case of (\ref{gasymptotics}).  
To derive (\ref{FTasympinfinity}) from
(\ref{gasymptotics}), we just have to conjugate by $\Om^{-1}$. 

We justify in this way the asymptotic formulae of Bourdeau in \cite{Bo95a} for the $\C P^n$ sigma-model (which corresponds to the solution $m=-x_0$ and the representation $\th=\la_{n+1}$).
For the solution $m=-x_0$ and the representation $\th=\wedge^k \la_{n+1}$ we obtain the asymptotic formulae appearing in \cite{Bo95a} for the $Gr_k(\C^{n+1})$ sigma-model. 

The soliton polytope for $Gr_k(\C^{n+1})$ is
the polytope in $\h_\sharp^\ast$ spanned by the weights
$x_{i_1}+\cdots+x_{i_k}$ where $0\le i_1<\cdots<i_k=n$. The vacua are the weight vectors $e_{i_1}\wedge\cdots\wedge e_{i_k}$. The projections of the weights to the Coxeter Plane are the complex numbers $-(\om^{i_1}+\dots+\om^{i_k})$.
The difference of two weights is a root if and only if their index sets have precisely $k-1$ common elements. When this happens, we have a soliton between the two vacua, whose particle type and mass are those associated to that root.  

We illustrate this  using the example of $Gr_3(\C^6)$, which was described somewhat mysteriously in section 4 of \cite{Bo95a}. 
\begin{figure}[h]
\begin{center}
\includegraphics[angle=90,origin=c,scale=0.4, trim=50 180 50 200]{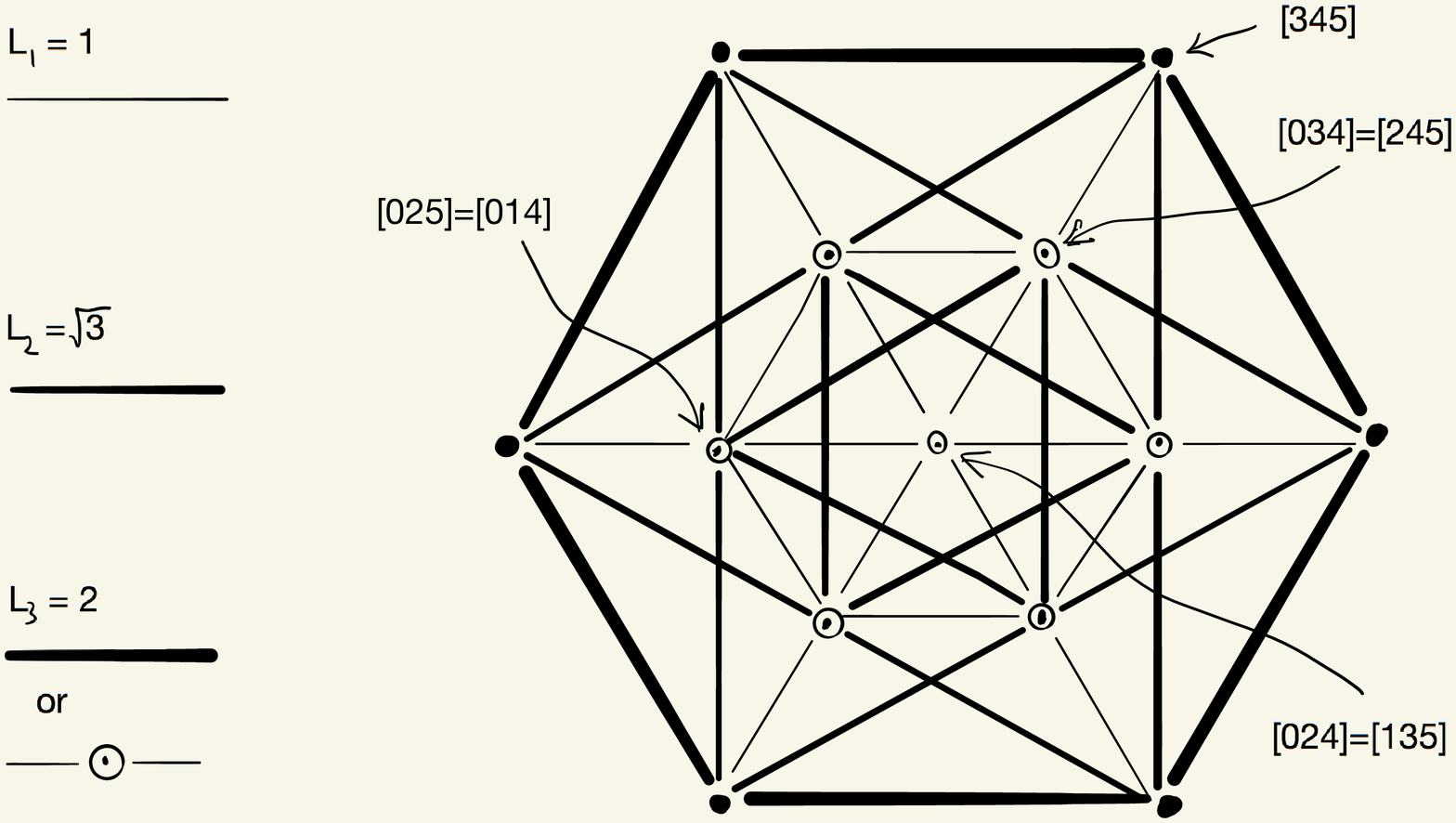}
\end{center}
\caption{The soliton structure of the $Gr_3(\C^6)$ sigma-model.}\label{gr36}
\end{figure}
For $\g=\sl_6\C$ there are $[\frac12(n+1)]=3$ particles, and
the particle data is as follows:
\begin{center}
$
\begin{array}{c|c|c|c}
\vphantom{\dfrac12}
 \text{particle} &  [x_0-x_1] & [x_0-x_2] &  [x_0-x_3]
\\
\hline
\vphantom{\dfrac12}
\text{mass}  & 2\sin\tfrac\pi6=1  & 2\sin\tfrac{2\pi}6=\sqrt 3  & 2\sin\tfrac{3\pi}6=2
 \\
\vphantom{\dfrac12}
\text{multiplicity} & \tbinom{6}{1}=6  & \tbinom{6}{2}=15 & \tbinom{6}{3}=20
\end{array}
$
\end{center}
The vertices of the soliton polytope, when projected to (our formulation of) the Coxeter Plane, produce the points $(x_i+x_j+x_k)(-d_{n+1})=
-(\om^i+\om^j+\om^k)$, where $\om=e^{2\pi\i/3}$.

Let us denote the state vector $e_i\wedge e_j\wedge e_k$ in $\wedge^3\C^6$ by $\vert ijk\rangle$ and its projection to the Coxeter Plane by $[ijk]$. As shown in 
 Figure \ref{gr36} we obtain $13$ distinct points:  the origin
 $[024]=[135]$, six vertices of a regular hexagon on a circle of radius $1$ (the Coxeter orbit of $[034]=[245]$ or $[025]=[014]$), and  six vertices of a regular hexagon on a circle of radius $2$ (the Coxeter orbit of $[345]$).
The Coxeter element $(543210)$ acts on the diagram by rotation through $-\pi/3$. 

The soliton structure can be read off from Figure \ref{gr36}.  For example, there is a soliton between the states $\vert 025\rangle$ and 
$\vert 245\rangle$ of type $[x_0-x_2]$ and multiplicity $15$.  Its mass is the distance between 
the points $[025]$ and 
$[245]$, namely $\sqrt 3$.  Examples of pairs connected by no solitons are
$\vert 024\rangle$ and 
$\vert 135\rangle$, or $\vert 034\rangle$ and $\vert 245\rangle$.

{\em
\noindent
Department of Mathematics\newline
Faculty of Science and Engineering\newline
Waseda University\newline
3-4-1 Okubo, Shinjuku, Tokyo 169-8555\newline
JAPAN
}

\end{document}